\documentclass[11pt]{amsart}
\hoffset= - 0.75 in

\usepackage{amsfonts, amsmath, amssymb, amsthm}
\usepackage{latexsym}

\newtheorem{thm}{Theorem}

\newtheorem{question}{Question}

\newtheorem{conj}{Conjecture}

\newtheorem{lem}{Lemma}
\newtheorem{cor}{Corollary}

\newcommand{\ZZ}{\mathbb{Z}}
\newcommand{\RR}{\mathbb{R}}

\newcommand{\ol}[1]{\overline{#1}}
\newcommand{\oa}[1]{\overleftrightarrow{#1}}
\newcommand{\bin}[2]{
	\left(
		\begin{array}{@{}c@{}}
			#1  \\  #2
		\end{array}
	\right)		}

\begin{document}

\title{A General Notion of Visibility Graphs}
\author[Develin, Hartke, and Moulton]{Mike Develin$^*$,
Stephen Hartke, and David Petrie
Moulton}
\thanks{$^*$Corresponding author. \textit{Address:} Department of Mathematics, UC-Berkeley, 
Berkeley, CA 94720-3840}
\email{develin@post.harvard.edu, hartke@math.rutgers.edu, moulton@idaccr.org}
\date{\today}

\begin{abstract}

We define a natural class of graphs by generalizing prior notions of
visibility, allowing the representing regions and sightlines to be arbitrary.
We consider mainly the case of compact connected representing regions, proving
two results giving necessary properties of visibility graphs, and giving
some
examples of classes of graphs that can be so represented. Finally, we give
some
applications of the concept, and we provide potential avenues for future
research in the area.

\end{abstract}

\maketitle

\section{Introduction}

At the AMS/MAA Joint Meetings in San Antonio in January of 1999, Alice
Dean presented a paper on rectangle-visibility graphs. A graph is
defined to be a rectangle-visibility graph if there exists a set of
closed rectangles on $\RR^2$, with sides parallel to the coordinate axes and disjoint except possibly for 
overlapping boundaries, such that each vertex corresponds to a rectangle and 
two vertices are joined by an edge if and only if there
exists an unobstructed horizontal or vertical corridor of positive width
connecting the two corresponding rectangles \cite{Dean2}.
There are many other concepts of visibility graphs; two of
the most commonly used are the line-segment definition (\cite{O'Rourke}, \cite{Agarwal}) and the simple polygon definition
(\cite{Abello}, \cite{Srin}).

In this paper, we consider visibility graphs corresponding to an extended
notion of visibility. Namely, we define a \textbf{sightline} to be a line
segment between points in two different regions intersecting no other
region. A graph $G$ (all graphs we consider will be finite and simple) is
then defined to be a \textbf{visibility graph} if there exists a set of
disjoint nonempty connected regions in $\RR^2$, with each region
corresponding to a vertex of $G$, such that two vertices are connected by
an edge in $G$ if and only if there exists a sightline between the
corresponding regions.  Throughout this paper, we will identify vertices
with the regions representing them.

Note that the regions are required to be connected. If we do not make this
requirement, all connected graphs turn out to be visibility graphs by
means of a simple construction. Let $G$ be an arbitrary connected graph; order the vertices $v_1,\ldots,v_n$ so 
that for all $1\le k\le n$, the induced 
subgraph $G|_{v_1,\ldots,v_k}$ is connected. Start with a disk representing $v_1$. For each $v_i$ in ascending 
order of $i$, cut a disk out of each area labeled $v_j$ for all $j<i$ such that $v_j$ and $v_i$ are adjacent in 
$G$. Label the interior of this disk $v_i$. It is then straightforward to check that this yields a visibility 
representation of $G$, as illustrated in Figure~1.

By requiring the representing regions to be connected, we eliminate
trivial constructions of this nature. In addition, we can require our
representing regions to have other properties such as openness,
compactness, and convexity; each of these yields a different theory of
visibility graphs. We can restrict the theory further by choosing our
regions from a severely limited list of objects, such as rectangles, disks
or even points. Naturally, with each restriction, the number of
representable graphs decreases.

We will be concerned primarily with \textbf{compact visibility
graphs}, graphs that can be represented as visibility graphs using compact 
connected subsets of
$\RR^2$. (A subset $S\subset \RR^2$ is \textbf{compact} if and only if it
is closed and bounded; this is equivalent to requiring every infinite 
sequence
$x_1,x_2,\ldots$ of points in $S$ to have a convergent subsequence.) 
This restriction allows
us to give and prove a necessary property of graphs in this class, while
at the same time preserving enough flexibility to facilitate 
visibility representations for a wide range of graphs.

\section{Convex Compact Visibility Graphs}

We begin by considering the case in which our representing regions are not 
only compact, but also convex. We
define a graph to be a \textbf{convex compact visibility graph}, or if it can be realized as a visibility
graph in $\RR^2$ with all representing regions both compact and convex. 
One of our main theorems is the following, which gives a necessary 
condition on such graphs.

\begin{thm}\label{premain}
If $G$ is a convex compact visibility graph, then every edge of $G$ must
be either a bridge or part of a $K_3$.
\end{thm}

The following lemma will be useful in the proof of the theorem.

\begin{lem}\label{conhull}

Suppose we have mutually visible regions $A$ and $B$. If any other region
intersects a segment joining a point of $A$ and a point of $B$, then the
edge $AB$ is part of a $K_3$ in the corresponding visibility graph.

\end{lem}

\begin{proof}
Take $a, a^\prime \in A$ and $b, b^\prime\in B$ so that $\overline{ab}$ is
a sightline and $\ol{a^\prime b^\prime}$ intersects some other region.
We have two cases: either $\overline{ab^\prime}$ is a sightline or it is 
not. It suffices to consider only the second case; if 
$\overline{ab^\prime}$ is a sightline, then switching the labels of the 
regions $A$ and $B$ and replacing $(a,b,b^\prime)$ by $(b^\prime, a, 
a^\prime)$ puts us in the situation of the second case. 

Consider the straight line path $\phi: [0,1]\rightarrow B$ from $b$ to
$b^\prime$; as $B$ is convex, the image of this path lies entirely in $B$.
Because all representing regions are compact, for
any given region $C$ intersecting some sightline $\overline{a\phi(t)}$
there exists some minimal value $t_C$
with $\overline{a\phi (t_C)}$ intersecting $C$ before $B$. (Note that if one
such line intersects $C$ before $B$, all such lines that intersect $C$ 
must
do so before $B$, as $B$ and $C$ are convex.) Since $\overline{a\phi (1)}$ is
occluded, at least one such interloping region exists. Pick the region $C$
with $t_C$ minimal, so that $\overline{a\phi (t)}$ is unobstructed for
$t<t_C$; if more than one such region exists, pick the one closest to $a$
along the line $\overline{a\phi (t_C)}$. Call the first point of
obstruction $c$. 
Then $A$ is adjacent to $C$ via the sightline $\overline{ac}$.
Furthermore, each point on $\overline{cb}$ lies on $\overline{a\phi (t)}$ for some
$t<t_C$ and thus must be unobstructed by the definition of $t_C$.
Therefore, $C$ is also adjacent to $B$ via the sightline
$\overline{cb}$, and so the edge $AB$ is part of the triangle $ABC$.
\end{proof}

We now proceed to the proof of Theorem~\ref{premain}.

\begin{proof}[Proof of Theorem~\ref{premain}]

Let $A$ and $B$ be mutually visible regions in a convex compact visibility
representation of $G$. Because $A$ and $B$ are compact convex regions, we
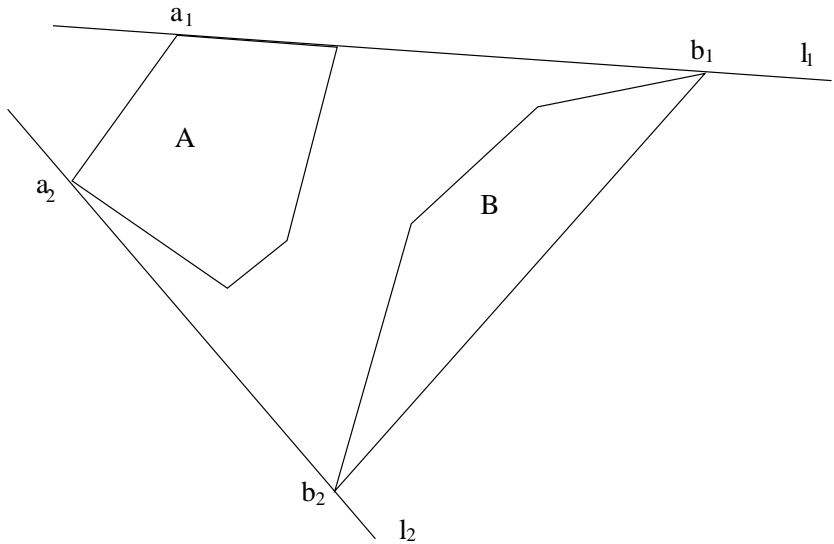
can find common external tangent lines $\ell_1$ and $\ell_2$ to $A$ and
$B$. Let $a_i$ (respectively, $b_i$) be a point in the intersection of $A$
(respectively, $B$) and $l_i$, as shown in Figure~2. If any other region
contains a point in the closed region bounded by $A$, $B$, and these two
tangent lines, then this point blocks some sightline between $A$ and $B$.
Therefore, we can apply Lemma~\ref{conhull} to conclude that the edge
corresponding to $AB$ is part of a triangle.

We may assume, therefore, that no such region exists. If every other
region is contained either in the area bounded by $A$, $\ell_1$, and
$\ell_2$ (and on the opposite side of $A$ from $B$) or that bounded by 
$B$, $\ell_1$, and $\ell_2$, it is clear that
$AB$ must be a bridge. Otherwise, we have one of two cases:

\textit{Case 1. Some region has nontrivial intersection with both
some $\ell_i$ and the open half-plane on the opposite side of that
line from $A$ and $B$.}

Without loss
of generality, we may assume that this line is $\ell_1$, and that the
region
intersects $\ell_1$ on the opposite side of $A$ from $B$.

Let $C$ be the first such region encountered by traveling along $\ell_1$
from $A$ away from $B$. Because $C$ extends at least to $\ell_1$, we can
construct, as in Figure~3, a sequence of points $\{c_i\}$ in $C$ so that
$\{c_i\}$ converges to $c\in C\cap \ell_1$, and each point $c_i$ is on the
opposite side of $\ell_1$ as $A$ and $B$.

Now, let $D$ be any region besides $C$ that intersects the open 
half-plane
on the other side of $l_1$ from $A$ and $B$. Then $D$ will not intersect
the segment $\ol{cb_1}$ by our choice of $C$. Since $D$ and $\ol{cb_1}$
are compact sets, the distance $\delta_D$ between them is positive; let
$\delta$ be the smallest such $\delta_D$. Then, if we pick $i$ large
enough so that $d(c_i,c)<\delta$, the sightlines $\ol{c_ia_1}$ and
$\ol{c_ib_1}$ must be unobstructed. Thus $ABC$ is a triangle.

\textit{Case 2. No such region exists; therefore, some region must lie
entirely on the opposite side of some $\ell_i$ from $A$ and $B$.}

We may suppose without loss of generality that this line is $\ell_1$.
Consider the
set of all regions $D_i$ lying on the opposite side of $\ell_1$ from $A$ and
$B$. Since each region is compact, and does not intersect $\ell_1$ itself,
each has some minimum distance $\delta_i>0$ to $\ell_1$. Choose a region $D$
with minimum distance $\delta$, and let $d\in D$ be such that
$d(d,\ell_1)=\delta$. Now, consider the sightlines $\overline{da_1}$ and
$\overline{db_1}$. These sightlines, which lie entirely on the opposite side
of $\ell_1$ from $A$ and $B$, must be unobstructed, for any intervening region
$E$ would necessarily have $d(E,\ell_1)<\delta$. Therefore $DAB$ is a
triangle.

Thus, in both cases, the edge $AB$ is part of a triangle. This completes the
proof of the theorem.
\end{proof}

The basic concept used in the proof of Theorem~\ref{premain} is the
division of $\RR^2$ into three areas: the convex hull of $A\cup B$, the
area blocked from $A$ by $B$ or from $B$ by $A$, and the remainder. Using
this division, we showed that if any region $C$ intersects either the
first area or the third, then $AB$ is part of a triangle, and if no such
region exists then $AB$ is a bridge. This division and the proof carry
through in the case where the representing regions are not required to be
convex. The details are similar to those of the proof of
Theorem~\ref{premain}, though slightly more technical, and are left to the
reader.

\begin{thm}\label{main}
If $G$ is a compact visibility graph, then every edge of $G$ is either a
bridge or part of a $K_3$.
\end{thm}

In light of this result, one might wonder whether or not the property of
convexity is irrelevant to the set of representable graphs; that is to 
say, whether all compact visibility 
graphs are compact convex visibility graphs. However, Figure~4 shows an
example of a graph that is a compact visibility graph, but not a convex 
compact visibility graph.

One might also wonder whether the converse of Theorem~\ref{main} is 
true, that
is, whether or not every graph satisfying the conclusion of 
Theorem~\ref{main} is
a compact visibility graph. This is easily shown to be false, however. To
construct a counterexample, we make use of the following lemma. 

\begin{lem}\label{ncon}
Suppose that $G$ is a visibility graph of any type, and $v$ is a vertex
all of whose neighbors are pairwise adjacent in $G$. Then $G-\{v\}$ is a
visibility graph of the same type. (By \textbf{type} of a visibility 
graph, we mean a type of restriction on the allowable regions in its 
representation.)
\end{lem}

\begin{proof}
Construct a visibility representation of $G$, and remove the region
corresponding to $v$. Consider any edge $AB$ in the resulting visibility
graph, with corresponding sightline $\ol{ab}$.  If we consider
$\overline{ab}$ in the original representation, the only region that can
possibly block it is the region corresponding to $v$. However, if this is
the case, then $A$ and $B$ are each connected to $v$ using subsegments of
$\overline{ab}$ as sightlines, so, by the hypothesis, $A$ and $B$ are
adjacent in the original graph. Furthermore, if $v$ does not block
$\overline{ab}$ then $A$ and $B$ are certainly adjacent in the original
graph. Thus the edge $AB$ also occurs in $G-\{v\}$. Clearly $G-\{v\}$ is
also a subgraph of the resulting visibility graph, so the two are
identical, and the lemma is proven. 
\end{proof}

Now, consider the graph $G$ shown in Figure~5, which satisfies the
conclusion of Theorem~\ref{main}. Suppose $G$ were a compact 
visibility
graph; by Lemma~\ref{ncon}, we could then remove the four outside vertices
to obtain a representation of $C_4$ using compact regions. However, no
such representation exists (as $C_4$ does not satisfy the conclusion of
Theorem~\ref{main}), and so $G$ cannot be a compact visibility graph. 

In the proof of Lemma~\ref{ncon}, we made explicit use of the geometry of
the visibility representation of $G$. Indeed, the geometry of this
representation yields a condition on the graph itself, as shown in the
following theorem.

\begin{thm}\label{k4}
Suppose that there exists a convex compact visibility representation of
$G$ such that two sightlines $\ol{ab}$ and $\ol{vw}$ intersect, with
$a,b,v,w$ contained in distinct regions $A,B,V,W$. Then the vertices $A$
and $B$ are part of a $K_4$ inside $G$. 
\end{thm}

\begin{proof}
If the regions $A$ and $B$ are both contained in the line $\oa{ab}$, or if
either region is a point, the proof is fairly straightforward and is left to
the reader. Suppose not; then, without loss of generality, there exists 
some point $a_0\in A$ not on $\oa{ab}$.

We define a point $a^\prime=(1-\alpha)a+\alpha a_0$, where $\alpha<1$ is
sufficiently small that $\overline{vw}$ still crosses $\overline{a^\prime
b}$ and $\ol{a^\prime b}$ is still a sightline. (If this latter condition
were impossible, we could construct a sequence of points in some other
region converging to $\overline{ab}$, and thus $\overline{ab}$ would also
not be a sightline.) We furthermore require that $B$ not be contained in
the line through $a^\prime$ and $b$.

Now, we have $\oa{a^\prime b}$ not tangent to $A$, and thus for $b^\prime$
sufficiently close to $b$ we have $\oa{a^\prime b^\prime}$ not tangent to
$A$. Pick $b_0\in B$ not on the line $\oa{a^\prime b}$.  Then, take
$b^\prime=(1-\beta)b+\beta b_0$, where $\beta$ is sufficiently small
that $\overline{a^\prime b^\prime}$ is a sightline and crosses $\ol{vw}$, 
and so that the
line $\oa{a^\prime b^\prime}$ is not tangent to $A$. Then $\oa{a^\prime
b^\prime}$ is tangent to neither $A$ nor $B$, for it has points from both
regions on either side of it.

Because of this construction, we can assume without loss of generality that
our extended sightline $\oa{ab}$ is tangent to neither $A$ nor $B$.

Let $\delta>0$ be the minimum distance between any pair of regions; these
distances are well-defined, as all regions are compact. For any line
segment $\overline{jk}$, let $f(\overline{jk})$ be the length of the
projection of $\ol{jk}$ onto a line perpendicular to $\ol{ab}$, and
$g(\overline{jk})$ be the length of the projection of $\ol{jk}$ onto a
line parallel to $\ol{ab}$; then $f$ and $g$ are always non-negative. Let
$x=\inf\{f(\ol{jk})\}$, where the infimum is taken over all sightlines
$\ol{jk}$ (between regions besides $A$ and $B$) intersecting $\ol{ab}$. It
is a simple consequence of compactness and the non-tangency of $\oa{ab}$
to $A$ and $B$ and compactness that we have $x>0$.

Now, since the regions are compact and the line $\oa{ab}$ is tangent to 
neither $A$ nor $B$, we can pick $\theta$ be sufficiently small so that
a line drawn from $a$ (respectively $b$) at an angle at most $\theta$ from
$\ol{ab}$ must intersect $B$ (respectively $A$) before any other region.
Pick $\varepsilon>0$ so that
$\varepsilon<\delta \sin \theta$. Let $\ol{cd}$ be a sightline crossing
$\ol{ab}$ with $f(\ol{cd})<x+\varepsilon$, and let $C$ and $D$ be the
regions containing $c$ and $d$ respectively.  We claim that $ABCD$ is a
$K_4$ in $G$. To verify this claim, we establish the existence of edges
$CA$, $CB$, $DA$, and $DB$.

Let $K$ be the triangle $cad$ in the plane, and consider the visibility graph $H$
defined by the regions $\{X\cap K\}$, where $X$ ranges over all the
regions intersecting $K$ from our original visibility representation of 
$G$; see Figure~6. As 
$K$ is
convex, and all of these new regions are contained in $K$, all 
relevant
sightlines will also be contained in $K$; furthermore, $H$ is clearly
a subgraph of $G$. Now, $H$ is a visibility graph, and $CD$ is an edge of
$H$. By Theorem~\ref{premain}, we know that $CD$ is either a bridge or 
part
of a triangle. 

Since $c$ and $d$ are on opposite sides of the unobstructed sightline
$\ol{ab}$, no region besides $A$ can intersect both $\overline{ca}$ 
and $\ol{ad}$. Taking the sequence of regions encountered along
$\ol{ca}$ gives a path from $C$ to $A$ in the graph $H$, and similarly
taking the sequence of regions along $\ol{ad}$ gives a path from $A$ to
$D$. Composing these yields a path from $C$ to $D$ not including the edge
$CD$, since $D$ cannot appear in the first part (as it appears in the 
second), and $C$ cannot appear in the second. Consequently, removing the 
edge $CD$ cannot disconnect the graph $H$, so $CD$ cannot be a bridge in 
$H$.

Therefore, $CD$ must be part of a triangle in $H$, so there must exist some region $E$ with both $CE$
and $DE$ edges of $H$. Unless $E=A$, $E$ must lie either entirely above
the unobstructed sightline $\ol{ab}\cap K$ or entirely below it. However,
if $E$ lies above this line, the sightline $\ol{d^\prime e}$ corresponding
to the edge $DE$ crosses $\ol{ab}$. Furthermore, as $e\in K$ and
the distance from $e$ to $c$ is at least $\delta$, it follows from our choice of $\varepsilon$
that $f(\ol{d^\prime e})<x$, a contradiction. Similarly, $E$ cannot lie
below $\ol{ab}\cap K$, as then $f(\ol{c^\prime e})<x$. Consequently, $E$
must be the region $A$, so that both of the edges $CA$ and $DA$ occur in
$H$ and, therefore, in $G$. Interchanging $A$ and $B$ shows that the edges 
$CB$
and $DB$ must also occur in $G$, so that $ABCD$ forms a $K_4$ as desired. 
\end{proof}

Theorem~\ref{k4} rules out an additional class of graphs, and places
constraints on possible representations of many more graphs. Noting that if a
visibility representation has no crossing sightlines, the corresponding 
graph is planar, we obtain the following corollary.

\begin{cor}
Every convex compact visibility graphs either is planar or contains a 
$K_4$.
\end{cor}

An example of a graph that this corollary excludes from being a convex
compact visibility graph is shown in Figure~7.
 
\section{Examples of Compact Visibility Graphs}
Theorem~\ref{main} shows that many simple graphs are not compact visibility
graphs; in particular, $C_n$ is not a compact visibility graph for $n>3$.
Furthermore, it is easy to check that every nonempty compact visibility
graph must be connected. Nevertheless, we can construct large families of 
compact visibility graphs.

\begin{thm}
If $G$ has a plane drawing such that all of the internal faces are
triangles, then $G$ is a compact visibility graph.
\end{thm}

\begin{proof}

Assume first that $G$ has no cut vertices, that is, vertices $v$ such that
$G-\{v\}$ is disconnected. We start with a plane drawing of $G$, which we
may assume has convex boundary face, and construct from this drawing a
compact visibility representation of $G$. Consider any edge $vw$. We break
$vw$ at an arbitrary point in its interior, and interlace the two halves 
of $vw$ together at this point, as shown in Figure~8.
We then take as the
representing region for $v$ the union of all the half-edges containing
$v$. It is clear that every edge $vw$ corresponds to
a sightline; furthermore, every sightline lies entirely in the interior of
some face of the drawing, and thus only connects regions corresponding to
the vertices of this face. As these vertices form a $K_3$ in the graph,
there are no extraneous edges added, and so this assignment of regions
shows that $G$ is in fact representable as a visibility graph. 

In the case where $G$ has a cut vertex $v$, we can split the vertex as in
Figure~9 and then represent the resulting graph as a
visibility graph as before, except that we do not do the above for the 
edge $v_1 v_2$. This yields a representation of $G$ as a visibility graph.

\end{proof}

A simple corollary of Theorem 4, as every tree has a plane drawing with no
internal faces, is the following.

\begin{cor}
All trees can be represented as compact visibility graphs.
\end{cor}

In fact, all trees are convex compact visibility graphs, or even compact
disk visibility graphs.

Since $C_n$ is not a compact visibility graph for $n>3$, one might suspect
that a planar graph $G$ is a visibility graph if and only if it has a
plane drawing such that all of the internal faces are triangles. This,
however, turns out to be false. Consider the case of $K_{1,1,n}$. This is
easily seen to be a compact visibility graph (as in Figure~10), but it is
a planar graph with no such drawing for $n\ge 3$.  The key point here is
that the induced subgraph on the four-sided face is not $C_4$; given this,
we conjecture that the following modification of the converse holds.

\begin{conj}
A planar graph $G$ is a compact visibility graph if and only if it has a
plane drawing such that for all internal faces $F$ of $G$, the 
subgraph $G_F$ induced by the vertices of $F$ is a compact visibility 
graph.
\end{conj}

\section{The non-compact case}
Throughout the previous two sections of this paper, we have assumed that
our
representing regions are compact. In this section, we give our rationale for
this assumption, as well as some comments about the non-compact theory.

Theorem~\ref{main} does not hold in the non-compact case. Indeed, the cycles
$C_n$ are visibility graphs, as shown by the two representations in 
Figure~11.

Furthermore, if we eliminate the hypothesis of compactness, we can even
realize disconnected graphs as visibility graphs. Let $A$ be the set
$\{(x,y)\,|\,0<x\le 1, y=\sqrt{x} \sin\frac{1}{x}\}$, $B$ be the set
$\{(x,x+y)\,|\,(x,y)\in A\}$, and $C$ be the single point $(0,0)$, as
represented in Figure~12. Then the only unblocked sightlines are between
$A$ and $B$, and so $\{A,B,C\}$ is a visibility representation of a
disconnected graph.

Combining the technique of interlacing presented in Figure~12 with the
second method for realizing $C_n$ presented in Figure~11, we can realize 
any
connected planar graph as a visibility graph. Given a face $F$, for each vertex $v$ we 
construct in the interior of $F$ interlocking spirals $s_{v,F}$, each having infinitely many 
turns; one of these spirals is shown in Figure~13. The region 
corresponding to $v$ then consists 
of the union over all faces $F$ containing $v$ of the spirals $s_{v,F}$.

While it is not true that we can represent any graph as a visibility graph
(consider the case of two isolated vertices), the fact that all
connected planar graphs can be realized as visibility graphs leads to the
following conjecture:

\begin{conj}\label{all}
Any connected graph is representable as a visibility graph.
\end{conj}

\section{Settings Other Than $\RR^2$}
It is natural to extend the concept of visibility graphs to 
higher-dimensional Euclidean
spaces. Specifically, one might consider the case of
$\RR^3$. Analogues to rectangle-visibility graphs exist in $\RR^3$; in
particular, Fekete and Meijer \cite{Fekete} consider the question of which
complete graphs are representable as box-visibility graphs. Alternatively, one
can require that the sightlines be strictly vertical, and that each
representing set lies in a horizontal plane; the question of such
representations has been considered by Alt, Godau, and Whitesides \cite{Alt}.

In $\RR^3$, Conjecture~\ref{all} is true (using our notion of visibility). 
Let $T$ be a
spanning tree for an arbitrary connected graph $G$; using the construction
shown in Figure~1 with spheres instead of circles, we first construct a
visibility representation of $T$. Then, to add an edge between two vertices
$V$ and $W$, we pick any points $v\in V$ and $w\in W$, and remove 
all points on the interior of the line segment $\ol{vw}$. It is clear that
this
will
only add
the single sightline from $V$ to $W$, as the line has thickness 0 and all
representing regions have positive thickness. Unlike in the
two-dimensional case, this construction does not disconnect any regions.
After all remaining
edges of $G$ have been added in this fashion, one is left with a visibility
representation of $G$.

In the compact case, we have results similar to those in $\RR^2$. With 
only slight modifications, the proof of Theorem~\ref{main} carries over
into $\RR^n$, so that Theorem~\ref{main} is true in general Euclidean 
spaces. This suggests the following two questions. 

\begin{question}
Is it true that all graphs representable as compact visibility graphs in
$\RR^3$ are also representable as compact visibility graphs in $\RR^2$?
\end{question}

\begin{question}
Does there exist a positive integer $n$ such that all graphs
representable as compact
visibility graphs in any Euclidean space are representable as compact
visibility graphs in $\RR^n$?
\end{question}

Two other venues to which visibility graphs generalize easily
are the torus $T^2$ and the sphere $S^2$; Mohar and Rosenstiehl 
\cite{Mohar}
briefly considered the former case with restricted sightlines. In order to
make sense of the general concept, we must define the notion of a straight
line in both settings. We define straight lines on the torus by considering it
as the quotient space $\RR^2/\ZZ^2$; on the sphere, we take the straight lines
to be the ordinary geodesics. Which graphs are representable as visibility
graphs in these new settings?

\section{Conclusion}
The generalized concept of a visibility graph lends itself to various
applications. Since the set of representing regions contains all of the
information about the graph itself, it can be used as an encoding of
information about the $n$ vertices of a graph $G$ and the $\bin{n}{2}$
connectedness relations
into a set of $n$ objects. If we require the representing regions to be
disks, we have this information using $3n$ coordinates: one triple
$(x,y,r)$ for each disk.

Visibility graphs also can be used naturally to represent various kinds of
networks. For instance, consider a network of army bases on a battlefield
that communicate with each other via torch signals (or, if you will,
lasers). The graph of this
network will then be the visibility graph represented by the layout of the
bases. Similarly, on a more global scale, one may wish to manipulate the
situation so that one can communicate in this fashion with one's allies and at
the same time intercept lines of communication between one's enemies. In 
this latter case of two rival alliances, it is appropriate to partition
the vertices into two classes, and only consider sightlines between two
regions of the same class; what pairs of networks can be realized in this
fashion? The concept of visibility graphs provides a natural
representation of graphs in $\RR^2$ that can be used to examine these
questions.

\section*{Acknowledgements}

The authors would like to thank Alice Dean for presenting the original talk which spurred the work
presented in this paper, and Joseph Gallian for helpful comments on a preliminary version of it; we
would also like to thank the referees for their suggestions. The first two authors were financially
supported by the University of Minnesota-Duluth and grants from the National Science Foundation (Grant
DMS-92820179) and the National Security Agency (Grant 904-00-1-0026), and by NSF Graduate Research
Fellowships. The third author was supported by the University of Wisconsin and the Center for
Communications Research in Princeton.

\end{document}